\documentclass{article}
\usepackage{amsmath,amsfonts,amssymb}
\usepackage[latin2]{inputenc}                                                   
\newtheorem{theorem}{Theorem}[section]
\newtheorem{definition}[theorem]{Definition}

\newtheorem{corollary}[theorem]{Corollary}

\newtheorem{example}[theorem]{Example}
\newtheorem{lemma}[theorem]{Lemma}

\newtheorem{remark}[theorem]{Remark}

\newcommand{\mC}{\mathbb C}

\newcommand{\mN}{\mathbb N}

\newcommand{\mZ}{\mathbb Z}

\newcommand{\mR}{\mathbb R}

\newcommand{\be}{\begin{eqnarray}}
\newcommand{\ee}{\end{eqnarray}}
\newcommand{\bd}{\begin{definition}}
\newcommand{\ed}{\end{definition}}

\newcommand{\br}{\begin{remark}}
\newcommand{\er}{\end{remark}}

\newcommand{\bt}{\begin{tabular}}
\newcommand{\et}{\end{tabular}}

\def\Sp{\mathop{\rm Sp}\nolimits}

\def\Mp{\mathop{\rm Mp}\nolimits}
\def\End{\mathop{\rm End}\nolimits}

\def\sp{\mathop{\mathfrak{sp}}\nolimits}
\def\sl{\mathop{\mathfrak{sl}}\nolimits}
\def\mp{\mathop{\mathfrak{mp}}\nolimits}

\def\Id{\mathop{\rm Id}\nolimits}

\input xypic
\input epsf
\input xy
\xyoption{all}

\begin{document}

\baselineskip13pt

\title{Symplectic twistor operator on $\mR^{2n}$ 
and the Segal-Shale-Weil representation}
\author{Marie Dost\'alov\'a, Petr Somberg
}
\date{}
\maketitle
\date{}
\abstract 
The aim of our article is the study of solution space of the symplectic twistor operator $T_s$ in symplectic spin geometry
on standard symplectic space $(\mR^{2n},\omega)$, which is the symplectic analogue of the twistor operator 
in (pseudo-)Riemannian spin geometry. In particular, we observe a substantial 
difference between the case $n=1$ of real dimension $2$ and the case of $\mR^{2n}$, $n>1$. For 
$n>1$, the solution space of $T_s$ is isomorphic to the Segal-Shale-Weil representation.   


\endabstract


\section{Introduction and Motivation}

In the case when the second Stiefel-Whitney class of a Riemannian manifold is trivial,
there is a double cover of the frame bundle and consequently there is an associated vector bundle
for the spinor representation of the spin structure group. There are two
basic first order invariant differential operators acting on spinor valued fields, namely the Dirac operator and
the twistor operator. Their spectral properties are reflected in the geometrical 
properties of the 
underlying manifold. In Riemannian geometry, the twistor equation appeared as an integrability 
condition for the
canonical almost complex structure on the twistor space, and it plays a prominent role 
in conformal differential geometry due to its
larger symmetry group. In physics, its solution space
defines infinitesimal isometries in Riemannian supergeometry. For an exposition with 
panorama of examples, cf. \cite{thom}, \cite{bfkg} and references therein.

The symplectic version of Dirac operator $D_s$ was introduced in \cite{KOS}, and its 
differential geometric properties were studied in \cite{crum}, \cite{MR2252919}, \cite{KAD}. 
The metaplectic Howe duality for $D_s$, introduced in \cite{bss}, allows to 
characterize the space of solutions for the symplectic Dirac operator $D_s$ on the 
(standard) symplectic space $(\mR^{2n},\omega)$. 

The aim of the present article is to study the symplectic twistor operator $T_s$
in the context of the metaplectic Howe duality, and consequently to determine its solution space 
on the standard symplectic space $(\mR^{2n},\omega)$. The operators $D_s$, $T_s$ were considered
from a different perspective in \cite{KAD}, \cite{Kr}, \cite{Kr1}. From an analytic point of view,
$T_s$ represents an overdetermined system of partial differential equations, acting 
on the space of polynomials valued in the vector space of the Segal-Shale-Weil
representation. From the point of view of representation theory, $T_s$ is $\mp(2n,\mR)$-invariant 
and the initial problem is dissolved by proper understanding of the interaction of 
$T_s$ with the generators $D_s,X_s$ of the Howe dual Lie algebra $\sl(2)$. 

As we shall see, as for $T_s$ there is a substantial difference between the situation for $n=1$ and 
$n>1$. Namely, there is in $\mathrm{Ker}(T_s)$ an infinite number of irreducible $\mp(2n,\mR)$-modules with 
different infinitesimal character for $n=1$, while for $n>1$ the kernel contains just
the Segal-Shale-Weil representation, a result of independent interest. This is the reason 
why we decided to treat the 
case $n=1$ in a separate paper \cite{ds} using different, more combinatorial approach,
which will be useful in complete understanding of the full infinite dimensional symmetry group of our operator.

The structure of our article goes as follows. In the first section, we review the subject of 
symplectic spin geometry and metaplectic Howe duality. In 
the second section, we start with the definition of the symplectic twistor operator $T_s$ and 
compute the space of polynomial solutions of $T_s$ on $(\mR^{2n},\omega)$. These results 
follow from the careful study of algebraic and differential consequences of $T_s$. 
In the last third section we indicate the collection of unsolved problems 
related to the topic of the present article.  

Throughout the article, we use the notation $\mN_0$ for the set of natural numbers including zero
and $\mN$ for the set of natural numbers without zero.

\subsection{Metaplectic Lie algebra $\mp(2n,\mR)$, symplectic Clifford algebra and a class of simple 
lowest weight modules for $\mp(2n,\mR)$}

In the present section we recall several algebraic and representation theoretical results used
in the next section 
for the analysis of the solution space of the symplectic twistor operator $T_s$, see e.g., 
\cite{bss}, \cite{crum}, \cite{fh}, \cite{MR2252919}, \cite{KAD}.

Let us consider $2n$-dimensional symplectic vector space 
$(\mR^{2n},\omega=\sum_{j=1}^{n}\epsilon^j\wedge \epsilon^{n+j}$), $n\in \mN$, 
and a symplectic basis $\{e_1,\ldots,e_{n},e_{n+1},\ldots, e_{2n}\}$ 
with respect to the non-degenerate two form $\omega\in\wedge^2({\mR^{2n}})^\star$.
Let $E_ {k,j}$ be the $2n\times 2n$ matrix with $1$ on the intersection of the 
$k$-th row and the $j$-th column and zero otherwise. 
The set of matrices 
$$X_{kj}=E_{k,j}-E_{n+j,n+k} \text{, } Y_{kj}=E_{k,n+j}+E_{j,n+k} \text{, } Z_{kj}=E_{n+k,j}+E_{n+j,k} ,$$
for $j,k= 1,\ldots, n$ is a basis of $\sp(2n,\mR)$, and
can be realized by first order differential operators
$$X_{kj}=x_j\partial_{x_k}-x_{n+k}\partial_{x_{n+j}} \text{, }\quad Y_{kj}=x_{n+j}\partial_{x_k}+x_{n+k}\partial_{x_j} \text{, }$$ $$Z_{kj}=x_j \partial_{x_{n+k}}+x_k \partial_{x_{n+j}}. $$

The metaplectic Lie algebra $\mp(2n,\mR)$ is the Lie algebra of the twofold group covering 
$\pi : \Mp(2n,\mR)\to \Sp(2n,\mR)$ of the symplectic Lie group $\Sp(2n,\mR)$. It can be realized 
by homogeneity two elements in 
the symplectic Clifford algebra $Cl_s(\mR^{2n},\omega)$, where the homomorphism $$\pi_\star : \mp(2n,\mR)\to \sp(2n,\mR)$$ 
is given by
\begin{eqnarray}
& & \pi_\star (e_k\cdot e_j)=-Y_{kj},
\nonumber \\
& & \pi_\star (e_{n+k}\cdot e_{n+j})=Z_{kj},
\nonumber \\
& & \pi_\star (e_{k}\cdot e_{n+j}+e_{n+j}\cdot e_k)=2X_{kj},
\end{eqnarray}
for $j,k=1,\ldots,n$.
\begin{definition}
The symplectic Clifford algebra $Cl_s(\mR^{2n},\omega)$ is an associative unital algebra over $\mC$,
given by the quotient of the tensor algebra $T(e_1,\ldots,e_{2n})$ by a
two-sided ideal $I\,\subset T(e_1,\ldots,e_{2n})$ generated by
$$
v_j\cdot v_k-v_k\cdot v_j=-i\omega (v_j,v_k)
$$ 
for all $v_j,v_k\in\mR^{2n}$, where $i\in\mC$ is the complex unit.
\end{definition}
The symplectic Clifford algebra $Cl_s(\mR^{2n},\omega)$ is isomorphic to the Weyl algebra $W_{2n}$ of complex 
valued algebraic 
differential operators on $\mR^n$, and the symplectic Lie algebra $\sp(2n,\mR)$ can be realized as 
a subalgebra of $W_{2n}$. In particular, the Weyl algebra is an associative algebra generated by 
$\{q_1,\ldots,q_{n},\partial_{q_1},\ldots,\partial_{q_n}\}$, the multiplication operator by $q_j$ and 
differentiation $\partial_{q_j}$, for $j=1,\ldots,n$, and the symplectic Lie algebra $\sp(2n,\mR)$ has a basis $\{-\frac{i}{2} q^2_j,
-\frac{i}{2} \frac{\partial^2}{\partial q^2_j}, q_j \frac{\partial}{\partial q_j}+\frac{1}{2}\}$, $j=1,\ldots,n$.

The symplectic spinor representation is an irreducible Segal-Shale-Weil representation of $Cl_s(\mR^{2n},\omega)$ on 
$L^2(\mR^n,e^{-\frac{1}{2}\sum_{j=1}^n q^2_j} dq_{\mR^n})$, the space of square integrable functions on 
$(\mR^n,e^{-\frac{1}{2}\sum_{j=1}^n q^2_j} dq_{\mR^n})$ with $dq_{\mR^n}$ the Lebesgue measure. 
Its action, the symplectic Clifford multiplication $c_s$, preserves the subspace
of $C^\infty$(smooth)-vectors given by the Schwartz space $S(\mR^n)$ of rapidly decreasing complex valued functions on
$\mR^n$ as a dense subspace. The space $S(\mR^n)$ can be regarded as a smooth (Frechet) globalization
of the space of $\tilde{K}=\widetilde{\mathrm{U}}(n)$-finite vectors in the representation, 
where $\tilde{K}\subset \Mp(2n,\mR)$ is the maximal compact subgroup given by the double cover of 
$K=\mathrm{U}(n)\subset \Sp(2n,\mR)$. Though we shall work in the smooth globalization $S(\mR^n)$, the representative
vectors are usually chosen to belong to the underlying Harish-Chandra module
of $\tilde{K}=\widetilde{\mathrm{U}}(n)$-finite vectors preserved by $c_s$. 

The function spaces associated to the Segal-Shale-Weil representation are supported on $\mR^n\,\subset\mR^{2n}$, a 
maximal isotropic subspace of $(\mR^{2n},\omega)$. 
In its restriction to $\mp(2n,\mR)$, $S(\mR^n)$ decomposes into two 
unitary representations realized on the subspace of even resp. odd functions:
\begin{eqnarray}\label{sshrepr}
\varrho :\mp(2n,\mR)\to \End(S(\mR^n)),
\end{eqnarray} 
where the basis vectors act by
\begin{eqnarray}
& & \varrho (e_j\cdot e_k)=iq_jq_k, 
\nonumber \\
& & \varrho (e_{n+j}\cdot e_{n+k})=-i\partial_{q_j}\partial_{q_k},
\nonumber \\
& & \varrho (e_j\cdot e_{n+j}+e_{n+j}\cdot e_j)=q_j\partial_{q_j}+\partial_{q_j} q_j.
\end{eqnarray}
for all $j,k=1,\dots ,n$.
In this representation the $Cl_s(\mR^{2n},\omega)$ acts 
on $L^2(\mR^n,e^{-\frac{1}{2}\sum_{j=1}^n q^2_j} dq_{\mR^n})$ by 
continuous unbounded operators with domain $S(\mR^{n})$. The space of
$\tilde{K}=\widetilde{\mathrm{U}}(n)$-finite vectors consists of even resp. odd 
homogeneity $\mp(2n,\mR)$-submodule 
$$
\{\mathrm{Pol}_{even}(q_1,\dots ,q_n)e^{-\frac{1}{2}\sum_{j=1}^n q^2_j}\},\quad
\{\mathrm{Pol}_{odd}(q_1,\dots ,q_n)e^{-\frac{1}{2}\sum_{j=1}^n q^2_j}\}.
$$
It is also an irreducible representation of 
$\mp(2n,\mR)\ltimes \mathfrak{h}(n)$, the semidirect product of $\mp(2n,\mR)$ and $(2n+1)$-dimensional 
Heisenberg Lie algebra $\mathfrak{h}(n)$ spanned by $\{e_1,\dots ,e_{2n},\Id\}$. In what follows we denote the Segal-Shale-Weil 
representation by 
${\mathcal S}$, and ${\mathcal S}\simeq {\mathcal S}_+\oplus{\mathcal S}_-$ as $\mp(2n,\mR)$-module.

Let us denote by $\mathrm{Pol}(\mR^{2n})$ the vector space of complex valued polynomials on $\mR^{2n}$, 
and by $\mathrm{Pol}_l(\mR^{2n})$ the subspace of homogeneity $l$ polynomials. The complex vector space
$\mathrm{Pol}_l(\mR^{2n})$ is as an irreducible $\mp(2n,\mR)$-module isomorphic to $\mathrm{S}^l(\mC^{2n})$, 
the $l$-th symmetric power of the complexification of the fundamental vector representation $\mR^{2n}$, $l\in\mN_0$.
 

\subsection{Segal-Shale-Weil
representation and the metaplectic Howe duality}

Let us review a representation-theoretical result of \cite{BL}, formulated in the 
opposite convention of highest weight metaplectic modules.
Let $\omega_1,\dots ,\omega_n$ be the fundamental weights of the Lie algebra $\sp(2n,\mR)$, 
and let $L(\lambda)$ denote the simple module over the universal
enveloping algebra ${\fam2 U}(\mp(2n,\mR))$ of $\mp(2n,\mR)$ generated by the highest weight vector of the weight 
$\lambda$.

Algebraically, the decomposition of the space of polynomial functions on $\mR^{2n}$ valued in the Segal-Shale-Weil 
representation corresponds to the tensor product 
of $L(-\frac{1}{2}\omega_n)$ resp. $L(\omega_{n-1}-\frac{3}{2}\omega_n)$ with 
symmetric powers $\mathrm{S}^k(\mC^{2n})$ of the fundamental vector 
representation $\mC^{2n}$ of $\sp(2n,\mR)$, $k\in\mN_0$. The following result is well known.

\begin{corollary}(\cite{BL})\label{decompositionssw}
We have for $L(-\frac{1}{2}\omega_n)$
\begin{enumerate}
\item 
In the even case $k=2l$ ($2l+1$ terms on the right-hand side): 
\begin{eqnarray*}
& & L(-\frac{1}{2}\omega_n)\otimes \mathrm{S}^k(\mC^{2n}) \simeq  
L(-\frac{1}{2}\omega_n) \oplus L(\omega_1+\omega_{n-1}-\frac{3}{2}\omega_n)
\nonumber \\
&  & \oplus L(2\omega_1-\frac{1}{2}\omega_n)\oplus L(3\omega_1+\omega_{n-1}-\frac{3}{2}\omega_n)
\oplus\dots
\nonumber \\
& & \oplus L((2l-1)\omega_1+\omega_{n-1}-\frac{3}{2}\omega_n)\oplus
L(2l\omega_1-\frac{1}{2}\omega_n),
\end{eqnarray*}

\item
In the odd case $k=2l+1$ ($2l+2$ terms on the right-hand side):
\begin{eqnarray*}
& & L(-\frac{1}{2}\omega_n)\otimes \mathrm{S}^k(\mC^{2n}) \simeq  
L(\omega_{n-1}-\frac{3}{2}\omega_n)\oplus L(\omega_1-\frac{1}{2}\omega_n)
\nonumber \\
& & \oplus L(2\omega_1+\omega_{n-1}-\frac{3}{2}\omega_n) 
\oplus L(3\omega_1-\frac{1}{2}\omega_n)
\oplus\dots
\nonumber \\
& & \oplus L(2l\omega_1+\omega_{n-1}-\frac{3}{2}\omega_n)
\oplus L((2l+1)\omega_1-\frac{1}{2}\omega_n),
\end{eqnarray*}
\end{enumerate}
We have for $L(\omega_{n-1}-\frac{3}{2}\omega_n)$
\begin{enumerate}
\item 
In the even case $k=2l$ ($2l+1$ terms on the right-hand side): 
\begin{eqnarray*}
& & L(\omega_{n-1}-\frac{3}{2}\omega_n)\otimes \mathrm{S}^k(\mC^{2n})  \simeq 
L(\omega_{n-1}-\frac{3}{2}\omega_n)\oplus L(\omega_1-\frac{1}{2}\omega_n)
\nonumber \\ 
&& \oplus L(2\omega_1+\omega_{n-1}-\frac{3}{2}\omega_n)
\oplus\dots\oplus L((2l-1)\omega_1-\frac{1}{2}\omega_n)
\nonumber \\
& & \oplus L(2l\omega_1+\omega_{n-1}-\frac{3}{2}\omega_n),
\end{eqnarray*}

\item
In the odd case $k=2l+1$ ($2l+2$ terms on the right-hand side):
\begin{eqnarray*}
& & L(\omega_{n-1}-\frac{3}{2}\omega_n)\otimes \mathrm{S}^k(\mC^{2n})  \simeq  
L(-\frac{1}{2}\omega_n)\oplus L(\omega_1+\omega_{n-1}-\frac{3}{2}\omega_n)
\oplus\dots
\nonumber \\
& & \oplus L(2l\omega_1-\frac{1}{2}\omega_n)\oplus
L((2l+1)\omega_1+\omega_{n-1}-\frac{3}{2}\omega_n).
\end{eqnarray*}
\end{enumerate}
\end{corollary}

A more geometrical reformulation of this statement is realized in the algebraic (polynomial) 
Weyl algebra and termed metaplectic Howe 
duality, \cite{bss}. 
The metaplectic analogue of the classical theorem on the separation of variables allows to decompose 
the space $\mathrm{Pol}(\mR^{2n})\otimes {\mathcal S}$ of complex polynomials valued in the Segal-Shale-Weil
representation under the action of $\mp(2n,\mR)$ into a direct sum of simple lowest weight 
$\mp(2n,\mR)$-modules
\begin{eqnarray}
\mathrm{Pol}(\mR^{2n})\otimes {\mathcal S}\simeq\bigoplus_{l=0}^\infty\bigoplus_{j=0}^\infty X_s^j{M}_l,
\end{eqnarray}
where we use the notation ${M}_l:={M}^+_l\oplus {M}^-_l$. This decomposition takes
the form of an infinite triangle 
\begin{eqnarray}\label{obrdekonposition}
\xymatrix@=11pt{\mathrm{P}_0 \otimes {\mathcal S} \ar@{=}[d] &  \mathrm{P}_1 \otimes {\mathcal S} \ar@{=}[d]& 
\mathrm{P}_2 \otimes {\mathcal S} \ar@{=}[d] & \mathrm{P}_3 \otimes {\mathcal S} \ar@{=}[d] & 
\mathrm{P}_4 \otimes {\mathcal S} \ar@{=}[d]& \mathrm{P}_5 \otimes {\mathcal S} \ar@{=}[d] \\
M_0 \ar[r] & X_s M_0 \ar @{} [d] |{\oplus} \ar[r] & X_s^2 M_0 \ar @{} [d] |{\oplus} \ar[r] & X_s^3 M_0 \ar @{} [d] |{\oplus}
 \ar[r] & X_s^4 M_0 \ar @{} [d] |{\oplus}\ar[r] & X_s^5 M_0 \ar @{} [d] |{\oplus} \\
& M_1 \ar[r] & X_s M_1 \ar @{} [d] |{\oplus}\ar[r] & X_s^2 M_1 \ar @{} [d] |{\oplus}
 \ar[r] & X_s^3 M_1 \ar @{} [d] |{\oplus}\ar[r] & X_s^4 M_1 \ar @{} [d] |{\oplus} \\
&& M_2 \ar[r] & X_s M_2 \ar @{} [d] |{\oplus}
 \ar[r] & X_s^2 M_2 \ar @{} [d] |{\oplus}\ar[r] & X_s^3 M_2 \ar @{} [d] |{\oplus} \\
&&& M_3 \ar[r] & X_s M_3 \ar @{} [d] |{\oplus}\ar[r] & X_s^2 M_3  \ar @{} [d] |{\oplus} \\
&&&& M_4 \ar[r] & X_s M_4 \ar @{} [d] |{\oplus}  \\
&&&&& M_5 
}
\end{eqnarray}

Let us now explain the notation used on the previous picture. First of all, we used the shorthand 
notation $\mathrm{P}_l=\mathrm{Pol}_l(\mR^{2n}),\, l\in\mN_0$, and all spaces and arrows on the picture have the following meaning.
Let $i\in\mC$ be the complex unit. The three operators
\begin{eqnarray}
& & X_s={} \sum_{j=1}^n (x_{n+j} \partial_{q_j} + i x_j q_j),
\nonumber \\
& & D_s={} \sum_{j=1}^n (i q_j \partial_{x_{n+j}} - \partial_{x_j}\partial_{q_j}),
\nonumber \\
& & E_s={} \sum_{j=1}^{2n} x_{j}\partial_{x_j},
\end{eqnarray}
where $D_s$ and $X_s$ acts on the previous picture horizontally but in the opposite 
direction, and fulfil the $\sl(2)$-commutation relations:
\begin{eqnarray}\label{comrel}
& & [E_s +n,D_s]=-D_s,
\nonumber \\
\label{slRels}
& & [E_s +n,X_s]=X_s,\\
& & [X_s,D_s]=i(E +n).\nonumber 
\end{eqnarray}
For the purposes of our article, we do not need the proper normalization of the generators $D_s,X_s,E_s$ making
the isomorphism with standard commutation relations in $\sl(2)$ explicit.

The elements of $\mathrm{Pol}(\mR^{2n})\otimes {\mathcal S}$ are called polynomial symplectic spinors.
Let $s\equiv s(x_1,\ldots,x_{2n},q_1,\ldots,q_n)\in \mathrm{Pol}(\mR^{2n})\otimes {\mathcal S}$, $h\in \Mp(2n,\mathbb{R})$ and 
$\pi (h)=g\in \Sp(2n,\mathbb{R})$ for the double covering map
$\pi: \Mp(2n,\mathbb{R})\rightarrow \Sp(2n,\mathbb{R})$. We define
the action of $\Mp(2n,\mR)$ on $\mathrm{Pol}(\mR^{2n})\otimes {\mathcal S}$ by
\begin{eqnarray}
\tilde{\varrho} (h)
s(x_1,\ldots,x_{2n},q_1,\ldots,q_n)=\varrho ( h) s(\pi(g^{-1}) (x_1,\ldots,x_{2n})^T,q_1,\ldots,q_n), 
\end{eqnarray} 
with $\varrho$ acting on the Segal-Shale-Weil representation via (\ref{sshrepr}).
Passing to the infinitesimal action, we get the operators representing the basis elements of $\mp(2n,\mR)$.
For example, we have for $j=1,\ldots,n$ 
\begin{align*}
\tilde{\varrho}(X_{jj})s= {}&\frac{d}{dt}\Big|_{t=0}\tilde{\varrho}(\exp(tX_{jj}))
s(x_1,\ldots,x_{2n},q_1,\ldots,q_n)\\
={} &\frac{d}{dt}\Big|_{t=0}
e^{\frac{t}{2}}
s(x_1,\ldots, x_j e^{-t},\ldots,x_{n+j} e^{t}\ldots,x_{2n},q_1,\ldots,q_je^{t},\ldots,q_n)\\
={} & \big(\frac{1}{2}- x_j \frac{\partial}{\partial x_j}+x_{n+j} \frac{\partial}{\partial x_{n+j}}+q_j \frac{\partial}{\partial q_j}\big) s(x_1,\ldots,x_{2n},q_1,\ldots,q_n).
\end{align*}
These operators satisfy the commutation relation of the Lie algebra $\mp(2n,\mathbb{R})$, 
and preserve the homogeneity in $x_1,\dots ,x_{2n}$. The operators $X_s$ and $D_s$ commute with operators 
$\tilde{\varrho}(X_{jk}),\tilde{\varrho}(Y_{jk})$ and $\tilde{\varrho}(Z_{jk})$, $j,k=1,\ldots,n$, hence
are $\mp(2n,\mR)$-intertwining differential operators. 

The action of $\mp(2n,\mR)\times \sl(2)$ generates the multiplicity free decomposition
of $\mathrm{Pol}(\mR^{2n})\otimes {\mathcal S}$ and the pair of Lie algebras in the product is called the metaplectic 
Howe dual pair. The operators $X_s$, $D_s$ acting on the previous picture horizontally
isomorphically identify the two neighboring $\mp(2n,\mR)$-modules. The modules $M_l$, $l\in\mN_0$,
on the most left diagonal of our picture are termed symplectic monogenics, and are characterized as 
$l$-homogeneous solutions of the symplectic Dirac 
operator $D_s$. Thus the decomposition is given as a vector space by tensor product of the symplectic monogenics
multiplied the by polynomial algebra of invariants $\mC[X_s]$. 


\section{The symplectic twistor operator $T_s$ and its solution space on $\mR^{2n}$}

We start with an abstract definition of the symplectic twistor operator $T_s$. 
Let $(M,\omega)$ be a 
$2n$-dimensional symplectic manifold, $p: P\to M$ 
a principal fiber $\Sp(2n,\mR)$-bundle of symplectic frames on $M$. 
A metaplectic structure
on $(M,\omega)$ is a principal fiber $\Mp(2n,\mR)$-bundle $\tilde{P}\to M$ 
together with bundle morphism $\tilde P\to P$, equivariant with respect to the
double covering $\Mp(2n,\mR)\to \Sp(2n,\mR)$. The manifold $(M,\omega)$ 
with a metaplectic structure is usually called symplectic spin manifold.
The symplectic manifold $M$
admits a metaplectic structure if and only if the second Stiefel-Whitney
class $w_2(M)$ is trivial, and the equivalence classes of metaplectic
structures are classified by $H^1(M,\mZ_2)$. There is a unique metaplectic structure on $(\mR^{2n},\omega)$.
\begin{definition}
Let $(M, \nabla ,\omega)$ be a symplectic spin manifold of dimension $2n$, $\nabla^s$ 
the associated symplectic spin covariant derivative
and $\omega\in C^\infty(M,\wedge^2T^\star M)$ a non-degenerate $2$-form such that $\nabla\omega=0$. We denote by 
$\{e_1,\ldots,e_{2n}\}$
a local 
symplectic frame. The symplectic twistor operator $T_s$ on $M$ is the first order differential operator 
$T_s$ acting on smooth symplectic spinors ${\mathcal S}$:
\begin{eqnarray}
& & \nabla^s :\, C^\infty(M,{\mathcal S})\longrightarrow T^\star M\otimes 
C^\infty(M,{\mathcal S}),\nonumber \\
& & T_s:=P_{\mathrm{Ker}(c)}\circ\omega^{-1}\circ\nabla^s:\, C^\infty(M,{\mathcal S})
\longrightarrow C^\infty(M,{\mathcal T}),
\end{eqnarray}
where ${\mathcal T}$ is the space of symplectic twistors, 
$T^\star M\otimes{\mathcal S}\simeq{\mathcal S}\oplus{\mathcal T}$, given by algebraic projection 
$$
P_{\mathrm{Ker}(c_s)}:T^\star M\otimes C^\infty(M,{\mathcal S})\longrightarrow C^\infty(M,{\mathcal T})
$$ 
on the kernel of the symplectic Clifford multiplication $c_s$. In the local symplectic coframe 
$\{\epsilon^1\}^{2n}_{j=1}$ dual to the symplectic frame $\{e_j\}^{2n}_{j=1}$ with respect to $\omega$,
we have the local formula for $T_s$:
\begin{align}\label{localformula}
T_s=\sum_{k=1}^{2n}\epsilon^k \otimes \nabla_{e_k}^s 
+ \frac{i}{n}\sum_{j,k,l=1}^{2n}\epsilon^l \otimes \omega^{kj} e_l \cdot e_j \cdot \nabla_{e_k}^s ,
\end{align} 
where $\cdot$ is the shorthand notation for the symplectic Clifford multiplication and $i\in\mC$ is the imaginary unit. 
We use the convention $\omega^{kj}=1$ for $j=k+n$ and $k=1,\dots ,n$, $\omega^{kj}=-1$ for $k=n+1,\dots ,2n$ 
and $j=k-n$, and $\omega^{kj}=0$ otherwise.
\end{definition}
The symplectic Dirac operator $D_s$ is defined as the image of the symplectic Clifford multiplication $c_s$,
and a symplectic spinor in the kernel of $D_s$ is called symplectic monogenic.
\begin{lemma}
The symplectic twistor operator $T_s$ is $\mp(2n,\mR)$-invariant. 
\end{lemma}
{\bf Proof:}
The property of invariance is a direct consequence of the equivariance of symplectic 
covariant derivative and invariance of algebraic projection $P_{\mathrm{Ker}(c_s)}$, and 
amounts to verify
\begin{equation}\label{invar}
T_s( \tilde{\varrho}(g)s)=\tilde{\varrho}(g)( T_s s)
\end{equation}
for any $g\in \mp(2n,\mR)$ and $s\in C^\infty(M,{\mathcal S})$.
Using the local formula (\ref{localformula}) for $T_s$ in a local chart 
$(x_1,\ldots,x_{2n})$, both sides of (\ref{invar}) are equal to
\begin{align*}
{}&\sum_{k=1}^{2n}\epsilon^k \otimes \varrho(g) 
\frac{\partial}{\partial x_k} \big[s\big(\pi(g)^{-1}x\big)\big]\\
&+ \frac{i}{n}\sum_{j,k,l=1}^{2n}\epsilon^l \otimes \omega^{kj} e_l\cdot  e_j\cdot  
\Big[\varrho(g) \frac{\partial}{\partial x_k} \big[s\big(\pi(g)^{-1}x\big)\big]\Big]
\end{align*}
and the proof follows.

\hfill
$\square$

In the case $M=(\mathbb{R}^{2n},\omega)$, the symplectic Dirac and the symplectic twistor operators are  
given by 
\begin{equation}
D_s=\sum_{j,k=1}^{2n}\omega^{kj}e_k\cdot \frac{\partial}{\partial x_j},
\end{equation}
\begin{equation}
T_s=\sum_{l=1}^{2n}\epsilon^l \otimes  \frac{\partial}{\partial x_l}  
+ \frac{i}{n}\sum_{j,k,l=1}^{2n}\epsilon^l \otimes \omega^{kj} e_l\cdot  e_j\cdot  \frac{\partial}{\partial x_k}
\label{TvR}
=\sum_{l=1}^{2n}\epsilon^l \otimes \Big( \frac{\partial}{\partial x_l}  
- \frac{i}{n} e_l\cdot D_s\Big), 
\end{equation}
and we restrict their action to the space of polynomial symplectic spinors. 

\begin{lemma}\label{prolong}
Let $s\in \mathrm{Pol}(\mathbb{R}^{2n},{\mathcal S})$ be a symplectic spinor in the solution space of 
the symplectic twistor operator $T_s$. Then $s$ is in the kernel of the square of the symplectic 
Dirac operator $D_s^2$.
\end{lemma}

{\bf Proof:} 
Let $s$ be a polynomial symplectic spinor in $\mathrm{Ker} (T_s)$,
\begin{equation}
T_s s=\sum_{l=1}^{2n}\epsilon^l \otimes \Big(\frac{\partial}{\partial x_l}  
- \frac{i}{n}  e_l\cdot D_s\Big)s=0,
\end{equation}
i.e.,
\begin{equation}\label{twistorsD}
\Big(\frac{\partial}{\partial x_l}  
- \frac{i}{n} e_l\cdot D_s\Big)s=0, \quad l=1,\dots, 2n.
\end{equation}
We apply to the last equation partial differentiation operator $\frac{\partial}{\partial x_m}$, multiply 
it by the skew symmetric form $\omega^{ml}$ and sum over $m=1,\dots ,2n$:
\begin{eqnarray}
\sum_{l,m=1}^{2n}\Big(\omega^{ml} \frac{\partial}{\partial x_m} \frac{\partial}{\partial x_l}  
- \frac{i}{n} \omega^{ml} e_l\cdot \frac{\partial}{\partial x_m} D_s\Big)s=0.
\end{eqnarray}
The first part is zero because of the skew-symmetry of $\omega$ and the symmetry in $m,l$, and the second part 
is (a non-zero multiple of) the square of the symplectic Dirac operator $D_s^2$. 
Hence 
\begin{eqnarray}
\sum_{l,m=1}^{2n}\frac{i}{n}\omega^{ml}  e_l\cdot  \frac{\partial}{\partial x_m} D_s s=-\frac{i}{n}D_s^2 s=0
\end{eqnarray}
and the proof is complete.

\hfill 
$\square$

Let us consider $\mp(2n,\mR)$-submodules in the split composition series
\begin{eqnarray}
\{0\}\subset \mathrm{Ker}(D_s)\subset \mathrm{Ker}(D_s^2)
\end{eqnarray}
with
\begin{eqnarray}
\mathrm{Ker}(D_s^2)\simeq \mathrm{Ker}(D_s)\oplus (\mathrm{Ker}(D_s^2)/ \mathrm{Ker}(D_s)),
\end{eqnarray}
and discuss which of them are in the solution space of the symplectic twistor operator $T_s$.  
We have 
\begin{eqnarray} 
\mathrm{Ker}(T_s)\simeq(\mathrm{Ker}(T_s)\cap \mathrm{Ker}(D_s))\oplus (\mathrm{Ker}(T_s)\cap (\mathrm{Ker}(D_s^2) / \mathrm{Ker}(D_s))).
\end{eqnarray}

\begin{lemma}
Let $n\in\mN$ and $s\in \mathrm{Pol}(\mathbb{R}^{2n}, {\mathcal S})$ be a symplectic spinor fulfilling
\begin{eqnarray}
s\in \mathrm{Ker}(T_s)\cap \mathrm{Ker}(D_s).
\end{eqnarray}
Then $s$ is a constant (i.e., independent of $x_1,\dots ,x_{2n}$) symplectic monogenic 
spinor. This 
is described by the following picture:
\begin{itemize}
\item
$\mathrm{Pol}(\mathbb{R}^{2n},{\mathcal S}_-)$:
\begin{eqnarray}
\xymatrix@=11pt{
*+[F]{M^-_0} \ar[rr] && {X_s M^-_0} \ar @{} [d] |{\oplus} \ar[rr] && X_s^2 M^-_0 \ar @{} [d] |{\oplus} \ar[rr] && X_s^3 M^-_0 \ar @{} [d] |{\oplus}
  \ar[rr] &&  \ldots\\
&& M^-_1 \ar[rr] && {X_s M^-_1} \ar @{} [d] |{\oplus}\ar[rr] && X_s^2 M^-_1 \ar @{} [d] |{\oplus}
  \ar[rr] && \ldots\\
&&& & M^-_2  \ar[rr] && {X_s M^-_2} \ar @{} [d] |{\oplus}
  \ar[rr] && \ldots\\
&&&&&& M^-_3 \ar[rr] &&\ldots \\
}
\end{eqnarray}

\item
$\mathrm{Pol}(\mathbb{R}^{2n},{\mathcal S}_+)$:
\begin{eqnarray}
\xymatrix@=11pt{
*+[F]{M^+_0} \ar[rr] && {X_s M^+_0} \ar @{} [d] |{\oplus} \ar[rr] && X_s^2 M^+_0 \ar @{} [d] |{\oplus} \ar[rr] && X_s^3 M^+_0 \ar @{} [d] |{\oplus}
  \ar[rr] &&  \ldots\\
&& M^+_1 \ar[rr] && {X_s M^+_1} \ar @{} [d] |{\oplus}\ar[rr] && X_s^2 M^+_1 \ar @{} [d] |{\oplus}
  \ar[rr] && \ldots\\
&&&& M^+_2  \ar[rr] && {X_s M^+_2} \ar @{} [d] |{\oplus}
  \ar[rr] && \ldots\\
&&&&&& M^+_3 \ar[rr] &&\ldots \\
}
\end{eqnarray}
\end{itemize}
\end{lemma}
{\bf Proof:}
Let $s\in \mathrm{Pol}(\mathbb{R}^{2n}, {\mathcal S})$ be a solution of the symplectic twistor operator, see ($\ref{twistorsD})$,
\begin{eqnarray}\nonumber
\Big(\frac{\partial}{\partial x_l}  
- \frac{i}{n} e_l\cdot D_s\Big)s=0, \quad l=1,\dots,2n,
\end{eqnarray}
and at the same time $s\in \mathrm{Ker}(D_s)$. This implies 
\begin{eqnarray}\nonumber
\frac{\partial}{\partial x_l}s=0, \quad l=1,\dots ,2n,
\end{eqnarray}
so $s$ is a constant symplectic spinor. The proof is complete.  

\hfill
$\square$

\begin{lemma}
Let $s\in \mathrm{Pol}(\mathbb{R}^{2n}, {\mathcal S})$ be a symplectic monogenic
spinor of homogeneity $h\in\mN_0$, i.e. $D_s(s)=0$. Then the symplectic spinor $X_s(s)$ 
has the following property:
\begin{enumerate}
\item
If $n=1$, then $X_s(s)$ is in the kernel of $T_s$ for any homogeneity $h\in\mN_0$. This 
is described by the following picture:
\begin{itemize}
\item
$\mathrm{Pol}(\mathbb{R}^2,{\mathcal S}_-)$:
\begin{eqnarray}
\xymatrix@=11pt{
{M^-_0} \ar[rr] && *+[F]{X_s M^-_0} \ar @{} [d] |{\oplus} \ar[rr] && X_s^2 M^-_0 \ar @{} [d] |{\oplus} \ar[rr] && X_s^3 M^-_0 \ar @{} [d] |{\oplus}
  \ar[rr] &&  \ldots\\
&& M^-_1 \ar[rr] && *+[F]{X_s M^-_1} \ar @{} [d] |{\oplus}\ar[rr] && X_s^2 M^-_1 \ar @{} [d] |{\oplus}
  \ar[rr] && \ldots\\
&&& & M^-_2  \ar[rr] && *+[F]{X_s M^-_2} \ar @{} [d] |{\oplus}
  \ar[rr] && \ldots\\
&&&&&& M^-_3 \ar[rr] &&\ldots \\
}
\end{eqnarray}

\item
$\mathrm{Pol}(\mathbb{R}^2,{\mathcal S}_+)$:
\begin{eqnarray}
\xymatrix@=11pt{
{M^+_0} \ar[rr] && *+[F]{X_s M^+_0} \ar @{} [d] |{\oplus} \ar[rr] && X_s^2 M^+_0 \ar @{} [d] |{\oplus} \ar[rr] && X_s^3 M^+_0 \ar @{} [d] |{\oplus}
  \ar[rr] &&  \ldots\\
&& M^+_1 \ar[rr] && *+[F]{X_s M^+_1} \ar @{} [d] |{\oplus}\ar[rr] && X_s^2 M^+_1 \ar @{} [d] |{\oplus}
  \ar[rr] && \ldots\\
&&&& M^+_2  \ar[rr] && *+[F]{X_s M^+_2} \ar @{} [d] |{\oplus}
  \ar[rr] && \ldots\\
&&&&&& M^+_3 \ar[rr] &&\ldots \\
}
\end{eqnarray}
\end{itemize}

\item

If $n>1$, then $X_s(s)$ is in the kernel of $T_s$ if and only if the homogeneity of $s$ is equal to $h=0$. This 
is described by the following picture:
\begin{itemize}
\item
$\mathrm{Pol}(\mathbb{R}^{2n},{\mathcal S}_-)$:
\begin{eqnarray}
\xymatrix@=11pt{
{M^-_0} \ar[rr] && *+[F]{X_s M^-_0} \ar @{} [d] |{\oplus} \ar[rr] && X_s^2 M^-_0 \ar @{} [d] |{\oplus} \ar[rr] && X_s^3 M^-_0 \ar @{} [d] |{\oplus}
  \ar[rr] &&  \ldots\\
&& M^-_1 \ar[rr] && {X_s M^-_1} \ar @{} [d] |{\oplus}\ar[rr] && X_s^2 M^-_1 \ar @{} [d] |{\oplus}
  \ar[rr] && \ldots\\
&&& & M^-_2  \ar[rr] && {X_s M^-_2} \ar @{} [d] |{\oplus}
  \ar[rr] && \ldots\\
&&&&&& M^-_3 \ar[rr] &&\ldots \\
}
\end{eqnarray}

\item
$\mathrm{Pol}(\mathbb{R}^{2n},{\mathcal S}_+)$:
\begin{eqnarray}
\xymatrix@=11pt{
{M^+_0} \ar[rr] && *+[F]{X_s M^+_0} \ar @{} [d] |{\oplus} \ar[rr] && X_s^2 M^+_0 \ar @{} [d] |{\oplus} \ar[rr] && X_s^3 M^+_0 \ar @{} [d] |{\oplus}
  \ar[rr] &&  \ldots\\
&& M^+_1 \ar[rr] && {X_s M^+_1} \ar @{} [d] |{\oplus}\ar[rr] && X_s^2 M^+_1 \ar @{} [d] |{\oplus}
  \ar[rr] && \ldots\\
&&&& M^+_2  \ar[rr] && {X_s M^+_2} \ar @{} [d] |{\oplus}
  \ar[rr] && \ldots\\
&&&&&& M^+_3 \ar[rr] &&\ldots \\
}
\end{eqnarray}
\end{itemize}
\end{enumerate} 
\end{lemma}
{\bf Proof:}
Let $s$ be a non-zero symplectic spinor in the kernel of $D_s$. The question is when
the system of partial differential equations acting on $s$,
\begin{eqnarray}
(\partial_{x_k} - \frac{i}{n} e_k\cdot D_s)X_s s=0,
\end{eqnarray}
holds for all $k=1,...,2n$. In other words, we ask when $X_s(s)$ is in the kernel
of the symplectic twistor operator. 
Let us multiply the $k$-th equation of the system by $x_k$ and sum over all $k$,
\begin{eqnarray}
(E_s - \frac{i}{n}X_sD_s)X_s s=0.
\end{eqnarray}
We use the $\sl(2)$-commutation relations for $X_s$, $D_s$ and for $E_s$, $X_s$, see (\ref{comrel}), and 
the fact that $s$ is in the kernel of $D_s$. This gives
\begin{eqnarray}
(E_sX_s-\frac{1}{n}X_sE_s-X_s)s=0.
\end{eqnarray}
Assuming that $s$ is of homogeneity $h$, $E_s s=h s$, the last equation reduces to 
\begin{eqnarray}
(h+1-\frac{h}{n}-1)X_s s=h(1-\frac{1}{n})X_s s=0.
\end{eqnarray}
Observe that $(1-\frac{1}{n})\not= 0$ for $n>1$, and $X_s$ is an $\mp(2n,\mR)$-intertwining map acting injectively
on $\mathrm{Pol}(\mathbb{R}^{2n},{\mathcal S})$ as a result of the metaplectic Howe duality (i.e., $s$ being non-zero 
implies $X_s(s)$ is non-zero.) Because $s$ is assumed to be non-zero, the last display implies that either 
\begin{enumerate}
\item
$h=0$ and $n\in\mN$ is arbitrary, or 
\item
$n=1$ and $h$ is arbitrary.
\end{enumerate}
A straightforward check for $n>1$ and the homogeneity $h=0$ gives
\begin{eqnarray}
(\partial_{x_k}-ie_k D_s)X_s s=\big(e_k +X_s\partial_{x_k}-\frac{i}{n}e_k E_s -e_k\big)s=0,
\end{eqnarray}
and in the case $n=1$ and arbitrary homogeneity we have
\begin{align}\label{odkaz33}
& (\partial_{x_1}-ie_1 D_s)X_s s
=\big(e_1+e_1 x_1\partial_{x_1}+e_2 x_2\partial_{x_1}-e_1 x_1 \partial_{x_1}-e_1 x_2\partial_{x_2}-e_1 \big)s
\nonumber \\
& =\big(x_2(e_2\partial_{x_1}-e_1 \partial_{x_2})\big)s=-x_2 D_s s=0.
\end{align}
As for the second component $(\partial_{x_2}-ie_2 D_s)$ of the symplectic twistor operator, the 
computation is analogous to the first one in (\ref{odkaz33}).
This completes the proof.

\hfill
$\square$

Let us summarize our results in the final theorem.
\begin{theorem}
The solution space of the symplectic twistor operator $T_s$ on standard symplectic space $(\mR^{2n},\omega)$
is given by $\mp(2n,\mR)$-modules in the boxes on the following pictures:
\begin{itemize}
\item
In the case $n=1$, we have for $\mathrm{Pol}(\mathbb{R}^2,{\mathcal S}_\pm)$:
\begin{eqnarray}
\xymatrix@=11pt{
*+[F]{M^\pm_0} \ar[rr] && *+[F]{X_s M^\pm_0} \ar @{} [d] |{\oplus} \ar[rr] && X_s^2 M^\pm_0 \ar @{} [d] |{\oplus} \ar[rr] && X_s^3 M^\pm_0 \ar @{} [d] |{\oplus}
  \ar[rr] &&  \ldots\\
&& M^\pm_1 \ar[rr] && *+[F]{X_s M^\pm_1} \ar @{} [d] |{\oplus}\ar[rr] && X_s^2 M^\pm_1 \ar @{} [d] |{\oplus}
  \ar[rr] && \ldots\\
&&& & M^\pm_2  \ar[rr] && *+[F]{X_s M^\pm_2} \ar @{} [d] |{\oplus}
  \ar[rr] && \ldots\\
&&&&&& M^\pm_3 \ar[rr] &&\ldots \\
}
\end{eqnarray}

\item
In the case $n>1$, we have for $\mathrm{Pol}(\mathbb{R}^{2n},{\mathcal S}_\pm)$:
\begin{eqnarray}
\xymatrix@=11pt{
*+[F]{M^\pm_0} \ar[rr] && *+[F]{X_s M^\pm_0} \ar @{} [d] |{\oplus} \ar[rr] && X_s^2 M^\pm_0 \ar @{} [d] |{\oplus} \ar[rr] && X_s^3 M^\pm_0 \ar @{} [d] |{\oplus}
  \ar[rr] &&  \ldots\\
&& M^\pm_1 \ar[rr] && {X_s M^\pm_1} \ar @{} [d] |{\oplus}\ar[rr] && X_s^2 M^\pm_1 \ar @{} [d] |{\oplus}
  \ar[rr] && \ldots\\
&&&& M^\pm_2  \ar[rr] && {X_s M^\pm_2} \ar @{} [d] |{\oplus}
  \ar[rr] && \ldots\\
&&&&&& M^\pm_3 \ar[rr] &&\ldots \\
}
\end{eqnarray}
\end{itemize}
\end{theorem}

An interested reader can easily verify the previous result for $n>1$ by 
taking a simple solution $s$ of $D_s$ of homogeneity at least one (it is sufficient 
to generate such a simple solution from dimension $n=1$ case) and check that
$X_s(s)\notin \mathrm{Ker}(T_s)$.
\begin{example}
In the case $n=2$ and the homogeneity $2$, the symplectic spinor 
\begin{eqnarray}
s=e^{-\frac{q_1^2 + q_2^2}{2}} (-i  x_1 x_2 +  x_1 x_4 +  x_2 x_3 + i  x_3 x_4)
\end{eqnarray}
is a solution of $D_s$. However, $X_s s$ is not a solution of the symplectic twistor 
operator $T_s$ because, for example, the first and the second components of $T_s X_s (s)$ are
nonzero: 
\begin{eqnarray*}
(T_ss)^1=\epsilon^1\otimes  e^{-\frac{q_1^2 + q_2^2}{2}} q_2 (x_2 + i x_4)^2\not= 0,\\
(T_ss)^2=\epsilon^2\otimes  e^{-\frac{q_1^2 + q_2^2}{2}} q_1 (x_1 + i x_3)^2\not= 0.
\end{eqnarray*}

\end{example}
It is much harder to verify the result 
$X_ss\in \mathrm{Ker}(T_s)$ for all polynomial symplectic spinors $s$, $s\in \mathrm{Ker}(D_s)$, in the case $n=1$, 
and we refer to \cite{ds} for a non-trivial combinatorial proof of 
this assertion. 
 
 We would like to emphasize that the kernel of our solution space realizes (for $n>1$) the 
Segal-Shale-Weil representation, a prominent $\Sp(2n,\mR)$-module with far-reaching impact on 
harmonic analysis.  
 

\section{Comments and open problems}

In the present section we comment on the results achieved in our 
article. 

First of all, notice that in the case of (both even and odd) orthogonal
algebras and the spinor representation as an orthogonal analogue of the Segal-Shale-Weil
representation, the solution space of the twistor operator for orthogonal Lie algebras on $\mR^n$ is 
given by two copies of the spinor representation, in complete analogy with the symplectic 
case, see \cite{bfkg} for $n\geq 3$. As for $n=2$, we were not able to find the required
result in the available literature, although we believe 
it is known to specialists. Here one half of the Dirac operator is the Dolbeault operator
and the twistor operator is its complex conjugate, while the opposite half of the Dirac and
twistor operators are their complex conjugates, respectively. The solution spaces for both halves of the 
twistor operator on $\mR^2$ are the complex linear spans of polynomials $\{z^j\}_{j\in\mN_0}$
and $\{\bar z^j\}_{j\in\mN_0}$, respectively, intersecting non-trivially in the constant polynomials. 
This is an orthogonal analogue of our results in symplectic
category, and indicates an infinite-dimensional symmetry group acting on the solutions
spaces of both symplectic Dirac and symplectic twistor 
operators in the real dimension $2$. 

Another observation is related to the proof of Lemma \ref{prolong} and its structure 
on curved symplectic manifolds. Let us consider 
a $2n$-dimensional metaplectic manifold $(M,\nabla,\omega)$, with $\nabla^s$ the 
metaplectic covariant derivative. Then a differential consequence of the symplectic 
twistor equation on $M$ is 
\begin{eqnarray}
\sum_{l,m=1}^{2n}\big(\omega^{ml}(\nabla^s_m,\nabla^s_l) - \frac{i}{n}D^2_s\big)s=0,
\end{eqnarray}
where the first term (skewing of the composition of metaplectic covariant derivatives) 
gives the action of the symplectic curvature of the symplectic
connection $\nabla^s$ on the space of sections of a metaplectic bundle on $M$. This 
equation should be thought of as a symplectic analogue of the equation 
\begin{eqnarray} 
D^2s=\frac{1}{4}\frac{n}{n-1}Rs,\quad n\geq 3
\end{eqnarray}   
in Riemannian spin geometry with $s$ a twistor spinor, $D$ the Dirac operator and
$R$ the scalar curvature of the Riemannian structure, cf. \cite{bfkg}. The prolongation 
of the symplectic twistor equation then constructs a linear connection and covariant 
derivative on the Segal-Shale-Weil representation, in such a way that the covariantly constant 
sections correspond to symplectic twistor spinors.  

Another ramification of our results is related to the higher order twistor equations
acting on symplectic spinors, with leading part
\begin{eqnarray}
s\mapsto \Pi^{hw}\nabla^s_{(i_1}\dots \nabla^s_{i_j)}s,
\end{eqnarray}
where $\Pi^{hw}$ is projection on the highest weight (or, Cartan) component and round 
brackets denote the symmetric part in the composition of symplectic covariant derivatives. For example,
their solution spaces can be studied by its interaction with the metaplectic Howe duality  
in an analogous way as we did in the first order case $(j=1)$.  


\vspace{0.2cm}
{\bf Acknowledgement:}  The authors gratefully acknowledge the support of the grant GA CR P201/12/G028 and SVV-2013-267317.

\vspace{0.3cm}

Marie Dost\'alov\'a, Petr Somberg

Mathematical Institute of Charles University,

Sokolovsk\'a 83, Praha 8 - Karl\'{\i}n, Czech Republic, 

E-mail: madost@seznam.cz, somberg@karlin.mff.cuni.cz.

\end{document}